\documentstyle[12pt,fleqn]{article}

\begin{document}

\title{Hyperbolic Complex Numbers in Two Dimensions}

\author{Silviu Olariu
\thanks{e-mail: olariu@ifin.nipne.ro}\\
Institute of Physics and Nuclear Engineering, Tandem Laboratory\\
76900 Magurele, P.O. Box MG-6, Bucharest, Romania}

\date{4 August 2000}

\maketitle

\abstract

A system of commutative hyperbolic complex numbers in 2 dimensions is studied
in this paper.  Exponential and trigonometric forms are obtained for the
hyperbolic twocomplex numbers.  Expressions are given for the elementary
functions of hyperbolic twocomplex variable.  The functions of a hyperbolic
twocomplex variable which are defined by power series are analytic.  Relations
of equality exist between partial derivatives of the real components a function
of a hyperbolic twocomplex variable.  The integral of a twocomplex function
between two points is independent of the path connecting the points.  A
hyperbolic twocomplex polynomial can be written as a product of linear or
quadratic factors, although the factorization may not be unique.

\endabstract

\section{Introduction}

A regular, two-dimensional complex number $x+iy$ 
can be represented geometrically by the modulus $\rho=(x^2+y^2)^{1/2}$ and 
by the polar angle $\theta=\arctan(y/x)$. The modulus $\rho$ is multiplicative
and the polar angle $\theta$ is additive upon the multiplication of ordinary 
complex numbers.

The quaternions of Hamilton are a system of hypercomplex numbers
defined in four dimensions, the
multiplication being a noncommutative operation, \cite{1} 
and many other hypercomplex systems are
possible, \cite{2a}-\cite{2b} but these interesting hypercomplex systems 
do not have all the required properties of regular, 
two-dimensional complex numbers which rendered possible the development of the 
theory of functions of a complex variable.

A system of hypercomplex numbers in 2 dimensions is described in this work,
for which the multiplication is associative and commutative, and for which
an exponential form and the concepts of analytic twocomplex
function and contour integration can be defined.
The twocomplex numbers introduced in this work have 
the form $u=x+\delta y$, the variables $x, y$ being real 
numbers.  
The multiplication rules for the complex units $1, \delta$ are
$1\cdot \delta=\delta, \delta^2=1$. In a geometric representation, the
twocomplex number $u$ is represented by the point $A$ of coordinates $(x,y).$
The product of two twocomplex numbers is equal to zero if both numbers are
equal to zero, or if one of
the twocomplex numbers lies on the line $x=y$ and the other on the 
line $x=-y$.

The exponential form of a twocomplex number, defined for $x+y>0,
x-y>0$, is $u=\rho\exp(\delta\lambda/2)$,
where the amplitude is $\rho=(x^2-y^2)^{1/2}$ and the argument is
$\lambda=\ln\tan\theta$, $\tan\theta=(x+y)/(x-y)$,
$0<\theta<\pi/2$. 
The trigonometric form of a twocomplex number is
$u=d\sqrt{\sin 2\theta}\exp\{(1/2)\delta \ln\tan\theta\}$,
where $d^2=x^2+y^2$.
The amplitude $\rho$ is equal to
zero on the lines $x=\pm y$. The division 
$1/(x+\delta y) $is possible provided that $\rho\not= 0$.
If $u_1=x_1+\delta y_1, 
u_2=x_2+\delta y_2$ are twocomplex
numbers of amplitudes and arguments $\rho_1,\lambda_1$ and respectively
$\rho_2, \lambda_2$, then the amplitude and the argument $\rho,
\lambda$ of
the product twocomplex number
$u_1u_2=x_1x_2+y_1y_2+\delta (x_1y_2+y_1x_2)$
are $\rho=\rho_1\rho_2,
\lambda=\lambda_1+\lambda_2$. 
Thus, the amplitude $\rho$ is a
multiplicative quantity
and the argument $\lambda$ is an additive quantity upon the
multiplication of twocomplex numbers, which reminds the properties of 
ordinary, two-dimensional complex numbers.

Expressions are given for the elementary functions of twocomplex variable.
Moreover, it is shown that the region of convergence of series of powers of
twocomplex variables is a rectangle  having the sides parallel to the bisectors
$x=\pm y$ . 

A function $f(u)$ of the twocomplex variable 
$u=x+\delta y$ can be defined by a corresponding 
power series. It will be shown that the function $f$ has a
derivative $\lim_{u\rightarrow u_0} [f(u)-f(u_0)]/(u-u_0)$ independent of the
direction of approach of $u$ to $u_0$. If the twocomplex function $f(u)$
of the twocomplex variable $u$ is written in terms of 
the real functions $P(x,y),Q(x,y)$ of real variables $x,y$ as
$f(u)=P(x,y)+\delta Q(x,y)$, then relations of equality 
exist between partial derivatives of the functions $P,Q$, and the 
functions $P,Q$ are solutions of the two-dimensional wave equation.

It will also be shown that the integral $\int_A^B f(u) du$ of a twocomplex
function between two points $A,B$ is independent of the 
path connecting the points $A,B$.

A polynomial $u^n+a_1 u^{n-1}+\cdots+a_{n-1} u +a_n $ can be
written as a 
product of linear or quadratic factors, although the factorization may not be
unique. 

This paper belongs to a series of studies on commutative complex numbers in $n$
dimensions.\cite{2c}
The twocomplex numbers described in this work are a particular case for 
$n=2$ of the polar hypercomplex numbers in $n$ dimensions.\cite{2c},\cite{2d}

\section{Operations with hyperbolic twocomplex numbers}

A hyperbolic complex number in two dimensions is determined by its two
components $(x,y)$. The sum 
of the hyperbolic twocomplex numbers $(x,y)$ and
$(x^\prime,y^\prime)$ is the hyperbolic twocomplex
number $(x+x^\prime,y+y^\prime)$. 
The product of the hyperbolic twocomplex numbers
$(x,y)$ and $(x^\prime,y^\prime)$ 
is defined in this work to be the hyperbolic twocomplex
number
$(xx^\prime+yy^\prime,
xy^\prime+yx^\prime)$.

Twocomplex numbers and their operations can be represented by  writing the
twocomplex number $(x,y)$ as  
$u=x+\delta y$, where $\delta$ 
is a basis for which the multiplication rules are 
\begin{equation}
1\cdot\delta=\delta,\: \delta^2=1 .
\label{1}
\end{equation}
Two twocomplex numbers $u=x+\delta y, 
u^\prime=x^\prime+\delta y^\prime$ are equal, 
$u=u^\prime$, if and only if $x=x^\prime, y=y^\prime. $
If $u=x+\delta y, 
u^\prime=x^\prime+\delta y^\prime$
are twocomplex numbers, 
the sum $u+u^\prime$ and the 
product $uu^\prime$ defined above can be obtained by applying the usual
algebraic rules to the sum 
$(x+\delta y)+ 
(x^\prime+\delta y^\prime)$
and to the product 
$(x+\delta y)
(x^\prime+\delta y^\prime)$,
and grouping of the resulting terms,
\begin{equation}
u+u^\prime=x+x^\prime+\delta(y+y^\prime),
\label{1a}
\end{equation}
\begin{equation}
uu^\prime=
xx^\prime+yy^\prime+
\delta(xy^\prime+yx^\prime)
\label{1b}
\end{equation}

If $u,u^\prime,u^{\prime\prime}$ are twocomplex numbers, the multiplication is
associative 
\begin{equation}
(uu^\prime)u^{\prime\prime}=u(u^\prime u^{\prime\prime})
\label{2}
\end{equation}
and commutative
\begin{equation}
u u^\prime=u^\prime u ,
\label{3}
\end{equation}
as can be checked through direct calculation.
The twocomplex zero is $0+\delta\cdot 0,$ 
denoted simply 0, 
and the twocomplex unity is $1+\delta\cdot 0,$ 
denoted simply 1.

The inverse of the twocomplex number 
$u=x+\delta y$ is a twocomplex number $u^\prime=x^\prime+\delta y^\prime$
having the property that
\begin{equation}
uu^\prime=1 .
\label{4}
\end{equation}
Written on components, the condition, Eq. (\ref{4}), is
\begin{equation}
\begin{array}{c}
xx^\prime+yy^\prime=1,\\
yx^\prime+xy^\prime=0.\\
\end{array}
\label{5}
\end{equation}
The system (\ref{5}) has the solution
\begin{equation}
x^\prime=\frac{x}{\nu} ,
\label{6a}
\end{equation}
\begin{equation}
y^\prime=
-\frac{y}{\nu} ,
\label{6b}
\end{equation}
provided that $\nu\not= 0$, where
\begin{equation}
\nu=x^2-y^2 .
\label{6e}
\end{equation}

The quantity $\nu$ can be written as
\begin{equation}
\nu=v_+v_-  ,
\label{7}
\end{equation}
where
\begin{equation}
v_+=x+y, \: v_-= x-y.
\label{8}
\end{equation}
The variables $v_+, v_-$ will be called canonical hyperbolic twocomplex
variables. 
Then a twocomplex number $u=x+\delta y$ has an inverse,
unless 
\begin{equation}
v_+=0 ,\:\:{\rm or}\:\: v_-=0.
\label{9}
\end{equation}

For arbitrary values of the variables $x,y$, the quantity $\nu$ can be
positive or negative. If $\nu\geq 0$, the quantity 
\begin{equation}
\rho=\nu^{1/2},\;\nu>0,
\label{9b}
\end{equation}
will be called amplitude of the twocomplex number $x+\delta y$.
The normals of the lines in Eq. (\ref{9}) are orthogonal to
each other. Because of conditions (\ref{9}) these lines will
be also called the nodal lines. 
It can be shown that 
if $uu^\prime=0$ then either $u=0$, or $u^\prime=0$, or 
one of the twocomplex numbers $u, u^\prime$ is of the form $x+\delta x$ and the
other is of the form $x-\delta x$.

\section{Geometric representation of hyperbolic twocomplex numbers}

The twocomplex number $x+\delta y$ can be represented by 
the point $A$ of coordinates $(x,y)$. 
If $O$ is the origin of the two-dimensional space $x,y$, the distance 
from $A$ to the origin $O$ can be taken as
\begin{equation}
d^2=x^2+y^2 .
\label{12}
\end{equation}
The distance $d$ will be called modulus of the twocomplex number $x+\delta y$. 

Since
\begin{equation}
(x+y)^2+(x-y)^2=2d^2 ,
\label{12b}
\end{equation}
$x+y$ and $x-y$ can be written as
\begin{equation}
x+y=\sqrt{2}d\sin\theta, \;x-y=\sqrt{2}d\cos\theta,
\label{12c}
\end{equation}
so that
\begin{equation}
x=d\sin(\theta+\pi/4),\; y=-d\cos(\theta+\pi/4) .
\label{12d}
\end{equation}

If $u=x+\delta y, u_1=x_1+\delta y_1,
u_2=x_2+\delta y_2$, and $u=u_1u_2$, and if
\begin{equation}
v_{j+}=x_j+y_j, \; v_{j-}= x_j-y_j , \;2d_j^2=v_{j+}^2+v_{j-}^2,\;
x_j+y_j=\sqrt{2}d_j\sin\theta_j, \;x_j-y_j=d_j\sqrt{2}\cos\theta_j,
\label{13}
\end{equation}
for $j=1,2$,
it can be shown that
\begin{equation}
v_+=v_{1+}v_{2+} ,\:\:
v_-=v_{1-}v_{2-}, \:\:\tan\theta=\tan\theta_1\tan\theta_2.
\label{14}
\end{equation}
The relations (\ref{14}) are a consequence of the identities
\begin{eqnarray}
(x_1x_2+y_1y_2)+(x_1y_2+y_1x_2)=(x_1+y_1)(x_2+y_2),
\label{15}
\end{eqnarray}
\begin{eqnarray}
(x_1x_2+y_1y_2)-(x_1y_2+y_1x_2)=(x_1-y_1)(x_2-y_2).
\label{16}
\end{eqnarray}

A consequence of Eqs. (\ref{14}) is that if $u=u_1u_2$, then
\begin{equation}
\nu=\nu_1\nu_2 ,
\label{19}
\end{equation}
where
\begin{equation}
\nu_j=v_{j+} v_{j-}, 
\label{20}
\end{equation}
for $j=1,2$.
If $\nu>0,\nu_1>0,\nu_2>0$, then
\begin{equation}
\rho=\rho_1\rho_2 ,
\label{20b}
\end{equation}
where
\begin{equation}
\rho_j=\nu_j^{1/2} , 
\label{20c}
\end{equation}
for $j=1,2$.

The twocomplex numbers
\begin{equation}
e_+=\frac{1+\delta}{2},\:
e_-=\frac{1-\delta}{2},\:
\label{21}
\end{equation}
are orthogonal,
\begin{equation}
e_+e_-=0,
\label{22}
\end{equation}
and have also the property that
\begin{equation}
e_+^2=e_+, \;e_-^2=e_-.
\label{23}
\end{equation}
The ensemble $e_+, e_-$ will be called the canonical hyperbolic twocomplex
base. 
The twocomplex number $u=x+\delta y$ can be written as
\begin{equation}
x+\delta y
=(x+y)e_++(x-y)e_-,  
\label{24}
\end{equation}
or, by using Eq. (\ref{8}),
\begin{equation}
u=v_+e_++v_- e_-,
\label{25}
\end{equation}
which will be called the canonical form of the hyperbolic twocomplex number.
Thus, if $u_j=v_{j+}e_++v_{j-} e_-, \:j=1,2$, and $u=u_1u_2$, then
the multiplication of the hyperbolic twocomplex numbers is expressed by the
relations (\ref{14}).

The relation (\ref{19}) for the product of twocomplex numbers can 
be demonstrated also by using a representation of the multiplication of the 
twocomplex numbers by matrices, in which the twocomplex number $u=x+\delta
y$ is represented by the matrix
\begin{equation}
\left(\begin{array}{cc}
x&y\\
y&x\\
\end{array}\right) .
\label{26}
\end{equation}
The product $u=x+\delta y$ of the twocomplex numbers
$u_1=x_1+\delta y_1, u_2=x_2+\delta y_2$, can be represented by the matrix
multiplication  
\begin{equation}
\left(\begin{array}{cc}
x&y\\
y&x\\
\end{array}\right) =
\left(\begin{array}{cc}
x_1&y_1\\
y_1&x_1\\
\end{array}\right) 
\left(\begin{array}{cc}
x_2&y_2\\
y_2&x_2\\
\end{array}\right) .
\label{27}
\end{equation}
It can be checked that
\begin{equation}
{\rm det}\left(\begin{array}{cccc}
x&y\\
y&x\\
\end{array}\right) =
\nu .
\label{28}
\end{equation}
The identity (\ref{19}) is then a consequence of the fact the determinant 
of the product of matrices is equal to the product of the determinants 
of the factor matrices. 

\section{Exponential and trigonometric forms of a twocomplex number}

The exponential function of the twocomplex variable $u$ can be defined by the
series
\begin{equation}
\exp u = 1+u+u^2/2!+u^3/3!+\cdots . 
\label{29}
\end{equation}
It can be checked by direct multiplication of the series that
\begin{equation}
\exp(u+u^\prime)=\exp u \cdot \exp u^\prime . 
\label{30}
\end{equation}
If $u=x+\delta y$, then  $\exp u$ can be calculated as
$\exp u=\exp x \cdot \exp (\delta y) $. According to Eq. (\ref{1}), 
\begin{equation}
\delta^{2m}=1, \delta^{2m+1}=\delta, 
\label{31}
\end{equation}
where $m$ is a natural number,
so that $\exp (\delta y)$ can be
written as 
\begin{equation}
\exp (\delta y) = \cosh y +\delta \sinh y .
\label{32}
\end{equation}
From Eq. (\ref{32}) it can be inferred that
\begin{eqnarray}
(\cosh t +\delta \sinh t)^m=\cosh mt +\delta \sinh mt .
\label{33}
\end{eqnarray}

The twocomplex numbers $u=x+\delta y$ for which
$v_+=x+y>0, \: v_-= x-y>0$ can be written in the form 
\begin{equation}
x+\delta y=e^{x_1+\delta y_1} .
\label{34}
\end{equation}
The expressions of $x_1, y_1$ as functions of 
$x, y$ can be obtained by
developing $e^{\delta y_1}$ with the aid of
Eq. (\ref{32}) and separating the hypercomplex components, 
\begin{equation}
x=e^{x_1}\cosh y_1 ,
\label{35}
\end{equation}
\begin{equation}
y=e^{x_1}\sinh y_1 ,
\label{36}
\end{equation}
It can be shown from Eqs. (\ref{35})-(\ref{36}) that
\begin{equation}
x_1=\frac{1}{2} \ln(v_+ v_- ) , \:
y_1=\frac{1}{2}\ln\frac{v_+}{v_- }.\:
\label{39}
\end{equation}
The twocomplex number $u$ can thus be written as
\begin{equation}
u=\rho\exp(\delta \lambda),
\label{40}
\end{equation}
where the amplitude is $\rho=(x^2-y^2)^{1/2}$ and the argument is
$\lambda=(1/2)\ln\{(x+y)/(x-y)\} $, for $x+y>0, x-y>0$.
The expression (\ref{40}) can be written with the aid of the variables $d,
\theta$, Eq. (\ref{12c}), as
\begin{equation}
u=\rho\exp\left(\frac{1}{2}\delta \ln\tan\theta\right),
\label{40b}
\end{equation}
which is the exponential form of the twocomplex number $u$, where
$0<\theta<\pi/2$. 

The relation between the amplitude $\rho$ and the distance $d$ is
\begin{equation}
\rho=d\sin^{1/2}2\theta.
\label{40c}
\end{equation}
Substituting this form of $\rho$ in Eq. (\ref{40b})
yields 
\begin{equation}
u=d\sin^{1/2}2\theta\exp\left(\frac{1}{2}\delta \ln\tan\theta\right),
\label{40d}
\end{equation}
which is the trigonometric form of the twocomplex number $u$.

\section{Elementary functions of a twocomplex variable}

The logarithm $u_1$ of the twocomplex number $u$, $u_1=\ln u$, can be defined
for $v_+>0, v_->0$ as the solution of the equation
\begin{equation}
u=e^{u_1} ,
\label{41}
\end{equation}
for $u_1$ as a function of $u$. From Eq. (\ref{40b}) it results that 
\begin{equation}
\ln u=\ln\rho+\frac{1}{2}\delta \ln\tan\theta
\label{42}
\end{equation}
It can be inferred from Eqs. (\ref{42}) and (\ref{14}) that
\begin{equation}
\ln(u_1u_2)=\ln u_1+\ln u_2 .
\label{43}
\end{equation}
The explicit form of Eq. (\ref{42}) is 
\begin{eqnarray}
\ln (x+\delta y)=
\frac{1}{2}(1+\delta)\ln(x+y)
+\frac{1}{2}(1-\delta)\ln(x-y),
\label{45}
\end{eqnarray}
so that the relation (\ref{42}) can be written with the aid of 
Eq. (\ref{21}) as 
\begin{equation}
\ln u = e_+\ln v_+ + e_-\ln v_-.
\label{44}
\end{equation}

The power function $u^n$ can be defined for $v_+>0, v_->0$ and real values
of $n$ as 
\begin{equation}
u^n=e^{n\ln u} .
\label{46}
\end{equation}
It can be inferred from Eqs. (\ref{46}) and (\ref{43}) that
\begin{equation}
(u_1u_2)^n=u_1^n\:u_2^n .
\label{47}
\end{equation}
Using the expression (\ref{44}) for $\ln u$ and the relations (\ref{22}) and
(\ref{23}) it can be shown that
\begin{eqnarray}
(x+\delta y)^n=
\frac{1}{2}(1+\delta)(x+y)^n
+\frac{1}{2}(1-\delta)(x-y)^n.
\label{48}
\end{eqnarray}
For integer $n$, the relation (\ref{48}) is valid for any $x,y$. The
relation (\ref{48}) for $n=-1$ is
\begin{equation}
\frac{1}{x+\delta y}=
\frac{1}{2}\left(\frac{1+\delta}{x+y}
+\frac{1-\delta}{x-y}
\right) .
\label{49}
\end{equation}

The trigonometric functions $\cos u$ and $\sin u $ of a twocomplex variable
are defined by the series
\begin{equation}
\cos u = 1 - u^2/2!+u^4/4!+\cdots, 
\label{50}
\end{equation}
\begin{equation}
\sin u=u-u^3/3!+u^5/5! +\cdots .
\label{51}
\end{equation}
It can be checked by series multiplication that the usual addition theorems
hold also for the twocomplex numbers $u_1, u_2$,
\begin{equation}
\cos(u_1+u_2)=\cos u_1\cos u_2 - \sin u_1 \sin u_2 ,
\label{52}
\end{equation}
\begin{equation}
\sin(u_1+u_2)=\sin u_1\cos u_2 + \cos u_1 \sin u_2 .
\label{53}
\end{equation}
The cosine and sine functions of the hypercomplex variables $\delta y$ can be
expressed as 
\begin{equation}
\cos\delta y=\cos y, \: \sin\delta y=\delta\sin y.
\label{54}
\end{equation}
The cosine and sine functions of a twocomplex number $x+\delta y$ can then be
expressed in terms of elementary functions with the aid of the addition
theorems Eqs. (\ref{52}), (\ref{53}) and of the expressions in  Eq. 
(\ref{54}). 

The hyperbolic functions $\cosh u$ and $\sinh u $ of the twocomplex variable
$u$ are defined by the series
\begin{equation}
\cosh u = 1 + u^2/2!+u^4/4!+\cdots, 
\label{57}
\end{equation}
\begin{equation}
\sinh u=u+u^3/3!+u^5/5! +\cdots .
\label{58}
\end{equation}
It can be checked by series multiplication that the usual addition theorems
hold also for the twocomplex numbers $u_1, u_2$,
\begin{equation}
\cosh(u_1+u_2)=\cosh u_1\cosh u_2 + \sinh u_1 \sinh u_2 ,
\label{59}
\end{equation}
\begin{equation}
\sinh(u_1+u_2)=\sinh u_1\cosh u_2 + \cosh u_1 \sinh u_2 .
\label{60}
\end{equation}
The $\cosh$ and $\sinh$ functions of the hypercomplex variables $\delta y$ can
be expressed as 
\begin{equation}
\cosh\delta y=\cosh y, \: \sinh\delta y=\delta\sinh y.
\label{61}
\end{equation}
The hyperbolic cosine and sine functions of a twocomplex number $x+\delta y$
can then be 
expressed in terms of elementary functions with the aid of the addition
theorems Eqs. (\ref{59}), (\ref{60}) and of the expressions in  Eq. 
(\ref{61}). 

\section{Twocomplex power series}

A twocomplex series is an infinite sum of the form
\begin{equation}
a_0+a_1+a_2+\cdots+a_n+\cdots , 
\label{64}
\end{equation}
where the coefficients $a_n$ are twocomplex numbers. The convergence of 
the series (\ref{64}) can be defined in terms of the convergence of its 2 real
components. The convergence of a twocomplex series can however be studied
using twocomplex variables. The main criterion for absolute convergence 
remains the comparison theorem, but this requires a number of inequalities
which will be discussed further.

The modulus of a twocomplex number $u=x+\delta y$ can be defined as 
\begin{equation}
|u|=(x^2+y^2)^{1/2} ,
\label{65}
\end{equation}
so that according to Eq. (\ref{12}) $d=|u|$. Since $|x|\leq |u|, |y|\leq |u|$,
a property of  
absolute convergence established via a comparison theorem based on the modulus
of the series (\ref{64}) will ensure the absolute convergence of each real
component of that series.

The modulus of the sum $u_1+u_2$ of the twocomplex numbers $u_1, u_2$ fulfils
the inequality
\begin{equation}
||u_1|-|u_2||\leq |u_1+u_2|\leq |u_1|+|u_2| .
\label{66}
\end{equation}
For the product the relation is 
\begin{equation}
|u_1u_2|\leq \sqrt{2}|u_1||u_2| ,
\label{67}
\end{equation}
which replaces the relation of equality extant for regular complex numbers.
The equality in Eq. (\ref{67}) takes place for $x_1=y_1,x_2=y_2$ or
$x_1=-y_1,x_2=-y_2$.
In particular
\begin{equation}
|u^2|\leq \sqrt{2}|u|^2 .
\label{68}
\end{equation}
The inequality in Eq. (\ref{67}) implies that
\begin{equation}
|u^m|\leq 2^{(m-1)/2}|u|^m .
\label{69}
\end{equation}
From Eqs. (\ref{67}) and (\ref{69}) it results that
\begin{equation}
|au^m|\leq 2^{m/2} |a| |u|^m .
\label{70}
\end{equation}

A power series of the twocomplex variable $u$ is a series of the form
\begin{equation}
a_0+a_1 u + a_2 u^2+\cdots +a_l u^l+\cdots .
\label{71}
\end{equation}
Since
\begin{equation}
\left|\sum_{l=0}^\infty a_l u^l\right| \leq  \sum_{l=0}^\infty
2^{l/2}|a_l| |u|^l ,
\label{72}
\end{equation}
a sufficient condition for the absolute convergence of this series is that
\begin{equation}
\lim_{l\rightarrow \infty}\frac{\sqrt{2}|a_{l+1}||u|}{|a_l|}<1 .
\label{73}
\end{equation}
Thus the series is absolutely convergent for 
\begin{equation}
|u|<c_0,
\label{74}
\end{equation}
where 
\begin{equation}
c_0=\lim_{l\rightarrow\infty} \frac{|a_l|}{\sqrt{2}|a_{l+1}|} ,
\label{75}
\end{equation}

The convergence of the series (\ref{71}) can be also studied with the aid of
the formula (\ref{48}) which, for integer values of $l$, is valid for any $x,
y, z, t$. If $a_l=a_{lx}+\delta a_{ly}$, and
\begin{eqnarray}
A_{l+}=a_{lx}+a_{ly}, A_{l-}= a_{lx}-a_{ly} , 
\label{76}
\end{eqnarray}
it can be shown with the aid of relations (\ref{22}) and (\ref{23}) that
\begin{equation}
a_l e_+=A_{l+} e_+, \: a_l e_-=A_{l-} e_-, \: 
\label{77}
\end{equation}
so that the expression of the series (\ref{71}) becomes
\begin{equation}
\sum_{l=0}^\infty \left(A_{l+}  v_+^l e_++
A_{l-}  v_-^le_-\right) ,
\label{78}
\end{equation}
where the quantities $v_+, v_-$
have been defined in Eq. (\ref{8}).
The sufficient conditions for the absolute convergence of the series 
in Eq. (\ref{78}) are that
\begin{equation}
\lim_{l\rightarrow \infty}\frac{|A_{l+1,+}||v_+|}{|A_{l+}|}<1,
\lim_{l\rightarrow \infty}\frac{|A_{l+1,-}||v_-|}{|A_{l-}|}<1.
\label{79}
\end{equation}
Thus the series in Eq. (\ref{78}) is absolutely convergent for 
\begin{equation}
|x+y|<c_+,\:
|x-y|<c_-,
\label{80}
\end{equation}
where 
\begin{equation}
c_+=\lim_{l\rightarrow\infty} \frac{|A_{l+}|}{|A_{l+1,+}|} ,\:
c_-=\lim_{l\rightarrow\infty} \frac{|A_{l-}|}{|A_{l+1,-}|} .
\label{81}
\end{equation}
The relations (\ref{80}) show that the region of convergence of the series
(\ref{78}) is a rectangle having the sides parallel to the bisectors $x=\pm y$.
It can be shown that $c_0=(1/\sqrt{2}){\rm min}(c,c^\prime)$, where
${\rm min}(c,c^\prime)$ designates the smallest of the numbers $c, c^\prime$.
Since $|u|^2=(v_+^2+v_-^2)/2$, it can be seen that the circular region of
convergence defined in Eqs. (\ref{74}), (\ref{75})
is included in the parallelogram defined in Eqs. (\ref{80}) and (\ref{81}).

\section{Analytic functions of twocomplex variables}

The derivative  
of a function $f(u)$ of the twocomplex variables $u$ is
defined as a function $f^\prime (u)$ having the property that
\begin{equation}
|f(u)-f(u_0)-f^\prime (u_0)(u-u_0)|\rightarrow 0 \:\:{\rm as} 
\:\:|u-u_0|\rightarrow 0 . 
\label{gs88}
\end{equation}
If the difference $u-u_0$ is not parallel to one of the nodal hypersurfaces,
the definition in Eq. (\ref{gs88}) can also 
be written as
\begin{equation}
f^\prime (u_0)=\lim_{u\rightarrow u_0}\frac{f(u)-f(u_0)}{u-u_0} .
\label{gs89}
\end{equation}
The derivative of the function $f(u)=u^m $, with $m$ an integer, 
is $f^\prime (u)=mu^{m-1}$, as can be seen by developing $u^m=[u_0+(u-u_0)]^m$
as
\begin{equation}
u^m=\sum_{p=0}^{m}\frac{m!}{p!(m-p)!}u_0^{m-p}(u-u_0)^p,
\label{gs90}
\end{equation}
and using the definition (\ref{gs88}).

If the function $f^\prime (u)$ defined in Eq. (\ref{gs88}) is independent of
the 
direction in space along which $u$ is approaching $u_0$, the function $f(u)$ 
is said to be analytic, analogously to the case of functions of regular complex
variables. \cite{3} 
The function $u^m$, with $m$ an integer, 
of the twocomplex variable $u$ is analytic, because the
difference $u^m-u_0^m$ is always proportional to $u-u_0$, as can be seen from
Eq. (\ref{gs90}). Then series of
integer powers of $u$ will also be analytic functions of the twocomplex
variable $u$, and this result holds in fact for any commutative algebra. 

If an analytic function is defined by a series around a certain point, for
example $u=0$, as
\begin{equation}
f(u)=\sum_{k=0}^\infty a_k u^k ,
\label{gs91a}
\end{equation}
an expansion of $f(u)$ around a different point $u_0$,
\begin{equation}
f(u)=\sum_{k=0}^\infty c_k (u-u_0)^k ,
\label{gs91aa}
\end{equation}
can be obtained by
substituting in Eq. (\ref{gs91a}) the expression of $u^k$ according to Eq.
(\ref{gs90}). Assuming that the series are absolutely convergent so that the
order of the terms can be modified and ordering the terms in the resulting
expression according to the increasing powers of $u-u_0$ yields
\begin{equation}
f(u)=\sum_{k,l=0}^\infty \frac{(k+l)!}{k!l!}a_{k+l} u_0^l (u-u_0)^k .
\label{gs91b}
\end{equation}
Since the derivative of order $k$ at $u=u_0$ of the function $f(u)$ , Eq.
(\ref{gs91a}), is 
\begin{equation}
f^{(k)}(u_0)=\sum_{l=0}^\infty \frac{(k+l)!}{l!}a_{k+l} u_0^l ,
\label{gs91c}
\end{equation}
the expansion of $f(u)$ around $u=u_0$, Eq. (\ref{gs91b}), becomes
\begin{equation}
f(u)=\sum_{k=0}^\infty \frac{1}{k!} f^{(k)}(u_0)(u-u_0)^k ,
\label{gs91d}
\end{equation}
which has the same form as the series expansion of 2-dimensional complex
functions. 
The relation (\ref{gs91d}) shows that the coefficients in the series expansion,
Eq. (\ref{gs91aa}), are
\begin{equation}
c_k=\frac{1}{k!}f^{(k)}(u_0) .
\label{gs92}
\end{equation}

The rules for obtaining the derivatives and the integrals of the basic
functions can 
be obtained from the series of definitions and, as long as these series
expansions have the same form as the corresponding series for the
2-dimensional complex functions, the rules of derivation and integration remain
unchanged.

If the twocomplex function $f(u)$
of the twocomplex variable $u$ is written in terms of 
the real functions $P(x,y),Q(x,y)$ of real
variables $x,y$ as 
\begin{equation}
f(u)=P(x,y)+\delta Q(x,y), 
\label{87}
\end{equation}
then relations of equality 
exist between partial derivatives of the functions $P,Q$. These relations
can be obtained by writing the derivative of the function $f$ as
\begin{eqnarray}
\lim_{\Delta x,\Delta y\rightarrow 0}\frac{1}{\Delta x+\delta \Delta y } 
\left[\frac{\partial P}{\partial x}\Delta x+
\frac{\partial P}{\partial y}\Delta y
+\delta\left(\frac{\partial Q}{\partial x}\Delta x+
\frac{\partial Q}{\partial y}\Delta y\right) 
\right] ,
\label{88}
\end{eqnarray}
where the difference $u-u_0$ in Eq. (\ref{gs89}) is 
$u-u_0=\Delta x+\delta\Delta y$. 
The relations between the partials derivatives of the functions $P, Q$ are
obtained by setting successively in   
Eq. (\ref{88}) $\Delta x\rightarrow 0, \Delta y=0$;
then $\Delta x= 0, \Delta y\rightarrow 0$. The
relations are 
\begin{equation}
\frac{\partial P}{\partial x} = \frac{\partial Q}{\partial y} ,
\label{89}
\end{equation}
\begin{equation}
\frac{\partial Q}{\partial x} = \frac{\partial P}{\partial y} .
\label{90}
\end{equation}

The relations (\ref{89})-(\ref{90}) are analogous to the Riemann relations
for the real and imaginary components of a complex function. It can be shown
from Eqs. (\ref{89})-(\ref{90}) that the components $P, Q$ are solutions
of the equations 
\begin{equation}
\frac{\partial^2 P}{\partial x^2}-\frac{\partial^2 P}{\partial y^2}=0,
\label{93}
\end{equation}
\begin{equation}
\frac{\partial^2 Q}{\partial x^2}-\frac{\partial^2 Q}{\partial y^2}=0,
\label{94}
\end{equation}
As can be seen from Eqs. (\ref{93})-(\ref{94}), the components $P, Q$ of
an analytic function of twocomplex variable are solutions of the wave 
equation with respect to the variables $x,y$.

\section{Integrals of twocomplex functions}

The singularities of twocomplex functions arise from terms of the form
$1/(u-u_0)^m$, with $m>0$. Functions containing such terms are singular not
only at $u=u_0$, but also at all points of the lines
passing through $u_0$ and which are parallel to the nodal lines. 

The integral of a twocomplex function between two points $A, B$ along a path
situated in a region free of singularities is independent of path, which means
that the integral of an analytic function along a loop situated in a region
free from singularities is zero,
\begin{equation}
\oint_\Gamma f(u) du = 0.
\label{105}
\end{equation}
Using the expression, Eq. (\ref{87})
for $f(u)$ and the fact that $du=dx+\delta  dy$, the
explicit form of the integral in Eq. (\ref{105}) is
\begin{eqnarray}
\oint _\Gamma f(u) du = \oint_\Gamma
[(Pdx+Qdy)+\delta(Qdx+Pdy)]
\label{106}
\end{eqnarray}
If the functions $P, Q$ are regular on the surface $\Sigma$
enclosed by the loop $\Gamma$,
the integral along the loop $\Gamma$ can be transformed with the aid of the
theorem of Stokes in an integral over the surface $\Sigma$ of terms of the form
$\partial P/\partial y -  \partial Q/\partial x$ and
$\partial P/\partial x -  \partial Q/\partial y$ 
which are equal to zero by Eqs. (\ref{89})-(\ref{90}), and this proves Eq.
(\ref{105}). 

The exponential form of the twocomplex numbers, Eq. (\ref{40}), contains no
cyclic variable, and therefore the concept of residue is not applicable to the
twocomplex numbers defined in Eqs. (\ref{1}).

\section{Factorization of twocomplex polynomials}

A polynomial of degree $m$ of the twocomplex variable 
$u=x+\delta y$ has the form
\begin{equation}
P_m(u)=u^m+a_1 u^{m-1}+\cdots+a_{m-1} u +a_m ,
\label{106b}
\end{equation}
where the constants are in general twocomplex numbers.
If $a_m=a_{mx}+\delta a_{my}$, and with the
notations of Eqs. (\ref{8}) and (\ref{76}) applied for $0, 1, \cdots, m$ , the
polynomial $P_m(u)$ can be written as 
\begin{eqnarray}
\lefteqn{P_m= \left[v_+^m 
+A_{1+} v_+^{m-1}+\cdots+A_{m-1,+} v_++ A_{m+} \right] e_+\nonumber}\\
&&+\left[v_-^m 
+A_{1-} v_-^{m-1} +\cdots+A_{m-1,-} v_-+ A_{m-}
\right]e_-. \nonumber\\
&&
\label{107}
\end{eqnarray}
Each of the polynomials of degree $m$ with real coefficients in Eq. (\ref{107})
can be written as a product
of linear or quadratic factors with real coefficients, or as a product of
linear factors which, if imaginary, appear always in complex conjugate pairs.
Using the latter form for the simplicity of notations, the polynomial $P_m$
can be written as
\begin{equation}
P_m=\prod_{l=1}^m (v_+-v_{l+})e_+
+\prod_{l=1}^m (v_--v_{l-})e_-,
\label{108}
\end{equation}
where the quantities $v_{l+}$ appear always in complex conjugate pairs, and the
same is true for the quantities $v_{l-}$.
Due to the properties in Eqs. (\ref{22}) and (\ref{23}),
the polynomial $P_m(u)$ can be written as a product of factors of
the form  
\begin{equation}
P_m(u)=\prod_{l=1}^m \left[(v_+-v_{l+})e_+
+(v_--v_{l-})e_-
\right].
\label{109}
\end{equation}
This relation can be written with the aid of Eqs. (\ref{25}) as
\begin{eqnarray}
P_m(u)=\prod_{l=1}^m (u-u_l),
\label{110}
\end{eqnarray}
where
\begin{eqnarray}
u_l=e_+ v_{l+} +e_- v_{l-}, 
\label{111}
\end{eqnarray}
for $l=1,...,m$.
The roots $v_{l+}$ and the roots $v_{l-}$ 
defined in Eq. (\ref{108}) may be ordered arbitrarily, which
means that Eq. (\ref{111}) gives sets of $m$ roots
$u_1,...,u_m$ of the polynomial $P_m(u)$, 
corresponding to the various ways in which the roots $v_{l+}, v_{l-}$
are ordered according to $l$ in each group. Thus, while the hypercomplex
components in Eq. (\ref{108}) taken 
separately have unique factorizations, the polynomial $P_m(u)$ can be written
in many different ways as a product of linear factors. 

If $P(u)=u^2-1$, the degree is $m=2$, the coefficients of the polynomial are
$a_1=0, a_2=-1$, the twocomplex components of $a_2$ are $a_{2x}=-1, a_{2y}=0$,
the components $A_{2+}, A_{2-}$ are $A_{2+}=-1, A_{2-}=-1$.
The expression, Eq. (\ref{107}), of $P(u)$ 
is $P(u)=e_+(v_+^2-1)+e_- (v_-^2-1)$, and
the factorization in Eq. (\ref{110}) is $u^2-1=(u-u_1)(u-u_2)$, where 
$u_1=\pm e_+\pm e_-, u_2=-u_1$. The factorizations are thus
$u^2-1=(u+1)(u-1)$ and $u^2-1=(u+\delta)(u-\delta)$. 
It can be checked that 
$(\pm e_+\pm e_-)^2=e_++e_-=1$.

\section{Representation of hyperbolic twocomplex complex numbers by irreducible
matrices} 

If the matrix in Eq. (\ref{26}) representing the twocomplex number $u$ is
called $U$, and 
\begin{equation}
T=\left(
\begin{array}{cc}
\frac{1}{\sqrt{2}}   &  \frac{1}{\sqrt{2}}   \\
-\frac{1}{\sqrt{2}}  &  \frac{1}{\sqrt{2}}   \\
\end{array}\right),
\label{112}
\end{equation}
it can be checked that 
\begin{equation}
TUT^{-1}=\left(
\begin{array}{cc}
x+y   & 0     \\
0     & x-y   \\
\end{array}\right).
\label{113}
\end{equation}
The relations for the variables $v_+=x+y, v_-=x-y$ for the
multiplication of twocomplex numbers have been written in Eq. (\ref{14}). The
matrix
$T U T^{-1}$  provides an irreducible representation
\cite{5} of the twocomplex numbers $u=x+\delta y$, in terms of matrices with
real coefficients.

\section{Conclusions}

An exponential form exists for the twocomplex numbers, involving the amplitude
$\rho$ and the argument $\lambda$. 
The twocomplex functions defined by series of powers are analytic, and
the partial derivatives of the components of the twocomplex functions are
closely related. The integrals of twocomplex functions are independent of path
in regions where the functions are regular. 
The polynomials of
tricomplex variables can be written as products of linear or quadratic factors.

\end{document}